\documentclass[10pt,b5paper,dvips]{amsart}
\usepackage{graphicx}
\usepackage[]{amsmath, amssymb, amsfonts, amsthm}

\newtheorem{prop}{Proposition}

\theoremstyle{definition}
\newtheorem{defin}{Definition}[section]
\newtheorem*{remark}{Remark}

\newcommand{\ra}{\rightarrow}
\newcommand{\wdt}{\widetilde}
\newcommand{\ol}{\overline}

\begin{document}
\title{A note on jet and geometric approach to higher order connections}
\author{Ma\"{\i}do Rahula}
\author{Petr Va\v{s}\'ik}


\begin{abstract}
We compare two ways of interpreting higher order connections.
The geometric approach lies in the decomposition of higher
order tangent space into the horizontal and vertical structures while the
jet--like approach considers a higher order connection as the
section of a jet prolongation of a fibered manifold. Particularly,
we use the Ehresmann prolongation of a general connection and study the result from the point of view of geometric theory. We pay attention to linear connections, too.
\end{abstract}\maketitle

\section{Introduction}

Several models of real objects are given as a smooth manifold and
one or more linear connections, e.g. material elasticity, see
\cite{e}. To obtain a manifold with just one characterization, one
has to consider a concept of a higher order connection. In this
paper, we recall the basic concepts of higher order connections from
both geometric and jet--like point of view, sections
\ref{jets} and \ref{con}. Let us note that the original ideas are those of C. Ehresmann, i.e. the definition of a connection by means of a horizontal distribution in a tangent space, the double fibered manifolds and holonomic and nonholonomic jets of fibered mappings. The first idea can be found in \cite{ehres}, the second one in \cite{ehres1}. The second idea was used for the case of vector bundles by Pradines, \cite{pradines}. Finally, the concept of holonomic and nonholonomic jets is widely studied in \cite{kms, kolmod, vasik, vasik1, virsik}. The first idea was extended in \cite{R3}, where the main formulae of higher order objects in multiple tangent spaces are derived, see also \cite{rahula}. In this paper we compare the jet--like and geometric approach.  We also recall a product of general
connections which leads to the so called Ehresmann prolongation and
show the reason why this operation is outstanding, especially
concerning semiholonomic connections, section \ref{ehres}. We study Ehresmann prolongation
of a connection from both points of view and show the analogues in both approaches.

\section{Jet prolongation of a fibered manifold}\label{jets}
Classical theory reads that $r$-th
nonholonomic prolongation $\wdt{J}^{r}Y$ of $Y$ is defined by the following
iteration:
\begin{enumerate}
\item $\wdt{J}^{1}Y=J^{1}Y,$ i.e. $\wdt J^1Y$ is a space of 1-jets of sections $M\ra Y$ over the target space $Y$.
\item $\wdt{J}^{r}Y=J^{1}(\wdt{J}^{r-1}Y\rightarrow M).$
\end{enumerate}
Clearly, we have an inclusion $ J^{r}Y\subset \wdt{J}^{r}Y$ given by
$j^{r}_{x}\gamma\mapsto j^{1}_{x}(j^{r-1}\gamma).$ \noindent
Further, $r$-th semiholonomic prolongation
$\overline{J}^{r}Y\subset\wdt{J}^{r}Y$ is defined by the following
induction. First, by $\beta_1=\beta_{Y}$ we denote the projection $J^1Y\ra Y$ and by $\beta_r=\beta_{\wdt{J}^{r-1}Y}$ the
projection $\wdt{J}^rY=J^1\wdt{J}^{r-1}Y\ra\wdt{J}^{r-1}Y,\ r=2,3,\ldots .$
If we set $\overline{J}^{1}Y=J^{1}Y $ and assume we have
$\overline{J}^{r-1}Y\subset\wdt{J}^{r-1}Y$ such that the restriction
of the projection $\beta_{r-1}:\wdt{J}^{r-1}Y\rightarrow
\wdt{J}^{r-2}Y$ maps $\overline{J}^{r-1}Y$ into
$\overline{J}^{r-2}Y,$  we can construct
$J^{1}\beta_{r-1}:J^{1}\overline{J}^{r-1}Y\rightarrow
J^{1}\overline{J}^{r-2}Y$ and define
$$\overline{J}^{r}Y=\{A\in J^{1}\overline{J}^{r-1}Y; \
\beta_{r}(A)=J^{1}\beta_{r-1}(A)\in\overline{J}^{r-1}Y\}.$$

If we denote by $\mathcal{FM}_{m,n}$ the category  with objects composed of fibered manifolds with $m$-dimensional bases
and $n$-dimensional fibres and morphisms formed by locally invertible fiber-preserving
mappings, then, obviously,
$J^r,\ol J^r$ and $\wdt J^r$ are bundle functors on $\mathcal{FM}_{m,n}$.

Alternatively, one can define the $r$-th order semiholonomic prolongation $\ol J^rY$ by means of natural target projections of nonholonomic jets, see \cite{virsik}.
For $r\geq q\geq 0$ let us denote by $\pi^r_{q}$ the target surjection $\pi^r_{q}:\wdt J^rY\ra \wdt J^qY$ with $\pi^r_r$ being the identity on $\wdt J^rY.$ We note
that the restriction of these projections to the subspace of semiholonomic jet prolongations will be denoted by the same symbol. By applying the functor
$J^k$ we have also the surjections $J^k\pi^{r-k}_{q-k}:\wdt J^r Y\ra\wdt J^q Y$ and, consequently, the element $X\in\wdt J^r Y$ is semiholonomic if and only if
\begin{equation}\label{def1}
(J^k\pi^{r-k}_{q-k})(X)=\pi^r_q(X)\  \text{for any integers}\ 1\leq k\leq q\leq r.
\end{equation}
In \cite{virsik}, the proof of this property can be found and the author finds it useful when handling semiholonomic connections and their prolongations.

Now let us recall local coordinates on higher order jet
prolongations of a fibered manifold.
$Y\ra M$. Let us denote by
$x^{i},\ i=1,\ldots,m$ the local coordinates on $M$ and $y^{p},\
p=1,\ldots,n$ the fiber coordinates of $Y\ra M$. We
recall that the induced coordinates on the holonomic prolongation
$J^{r}Y$ are given by $(x^{i},y^{p}_{\alpha}),$ where $\alpha$ is a
multiindex of range $m$ satisfying $|\alpha|\leq r.$ Clearly, the coordinates
$y^{p}_{\alpha}$ on $J^{r}Y$ are characterized by the complete
symmetry in the indices of $\alpha$. Having the nonholonomic
prolongation $\wdt{J}^{r}Y$ constructed by the iteration, we define the local coordinates inductively as follows:
\begin{itemize}
\item[1)]
Suppose that the induced coordinates on $\wdt{J}^{r-1}Y$ are of the form
$$(x^i,y^p_{k_1\ldots k_{r-1}}),\ k_1,\ldots, k_{r-1}=0,1,\ldots,m.$$
\item[2)]
We define the induced coordinates on $\wdt{J}^r Y$ by
$$(x^i,y^p_{k_1\ldots k_{r-1}0}=y^p_{k_1\ldots k_{r-1}},y^p_{k_1\ldots k_{r-1}i}=\frac{\partial}{\partial x^i}y^p_{k_1\ldots k_{r-1}}).$$
\end{itemize}
It remains to describe coordinates on the semiholonomic prolongation $\overline{J}^{r}Y$. Let
$(k_{1},\ldots,k_{r})$, $\ k_{1},\ldots,k_{r}
=0,1,\ldots,m$ be a sequence of indices and denote by $\langle
k_{1},\ldots,k_{s}\rangle$, $s\leq r$ the sequence of
non-zero indices in $(k_{1},\ldots,k_{r})$ respecting the order.
Then the definition of $\overline{J}^{r}Y$ reads that
the point $(x^{i},y^{p}_{k_{1},\ldots,k_{r}})\in \wdt{J}^{r}Y$
belongs to $\overline{J}^{r}Y$ if and only if
$y^{p}_{k_{1},\ldots,k_{r}}=y^{p}_{l_{1},\ldots,l_{r}}$ whenever
$\langle k_{1},\ldots,k_{r}\rangle = \langle
l_{1},\ldots,l_{r}\rangle$

\section{Iterated tangents}
Another concept, in this paper called geometric, of a connection rises from the theory of iterated tangent spaces. Let us recall that the bundle $T^kM\rightarrow T^{k-1}M$ is equipped with the structure of a
$k$-fold vector bundle. Particularly,
$T^kM$ admits $k$ different projections to
$T^{k-1}M$\,,\[ \,\,\rho_s:=
T^{k-s}\pi_s:T^{k}M\rightarrow T^{k-1}M,\,\,
\]
where $\pi_s$ is the natural projection $T^sM\rightarrow T^{s-1}M\,,
s=1,2,\ldots, k$. Each projection defines a vector bundle with basis
$T^{k-1}M$ and the total space is composed of
$2^{k-1}n$-dimensional vector spaces as fibers.
The local coordinates on the neighborhoods
\[
T^sU\subset T^{s}M, \quad\text{where}\quad
T^{s-1}U=\pi_s(T^sU),\,\,s=1,2,\ldots, k,
\]
are derived from coordinates, or coordinate mappings, $(u^i)$, which are given on the neighborhood
$U\subset M$\,:\\

$U$:\qquad $(u^i),\,\,i=1,2,\dots,\, n,$

$TU$:\,\,\quad $(u^i, u_1^i),$\quad where\quad $u^i:=
u^i\circ\pi_1,\,\,u^i_1:= du^i,$

$T^2U$:\quad $(u^i, u_1^i,u^i_2,u^i_{12}),$

 where\quad $u^i:=
u^i\circ\pi_1\pi_2,\,\,u^i_1:= du^i\circ\pi_2,\,\,\, u^i_2:=
d(u^i\circ\pi_1),\,\,\,u^i_{12}:= d^2u^i,$

etc.\\

\begin{prop} Coordinate mappings given on the
neighborhood $T^{s-1}U$ induce coordinate mappings
on the neighborhood $T^sU$ with respect to the projection
$\pi_s$ by adding the
differentials of these mappings.
\end{prop}
\noindent
Local coordinates are obtained by the following principle:\\
to the coordinates of a point of a manifold we attach the coordinates of the vector tangent to the manifold
at that point. We use the following notation\,:\,\,
the coordinates of a neighborhood $T^kU$ consist of two copies of local coordinates on
  $T^{k-1}U$ where the second copy is equipped with
  an additional subscript $k$\,.  This principle is suitable in the sense that the coordinates
  with index $s$
  are recognized as the fiber coordinates for projections
 $\rho_s,\,\, s=1,2,\ldots,k$, i.e. the coordinates with index
 $s$ disappear after the application of projection $\rho_s$.\\

The coordinate form of the three projections
$\rho_s:T^{3}U\rightarrow T^2U,\,\,\, s=1,2,3,$ is given by the
following diagram:

 \[
(u^i, u_1^i,u^i_2,u^i_{12},u^i_3, u_{13}^i,u^i_{23},u^i_{123})
\]
\[
_{\rho_1}\swarrow\qquad\quad_{\rho_2}\downarrow\qquad\qquad\quad\searrow\,_{\rho_3}
\]
\[
(u^i, u^i_2,u^i_3,u^i_{23})\qquad(u^i, u^i_1,u^i_3
,u^i_{13})\qquad(u^i, u^i_1,u^i_2,u^i_{12}).
 \]\\

 \begin{remark}
 Let us note that the semiholonomity condition is connected to the notion of the osculating bundle, see \cite{rahula}, and can be defined as the equalizer of all possible projections, which corresponds to \eqref{def1}.
 \end{remark}

\section{Connections}\label{con}

We start with the jet--like approach to connections. This rather structural description is quite suitable for determining natural operators on connections, for details see \cite{kms}.
\begin{defin}
A general connection on the fibered manifold $Y\ra M$ is a section $\Gamma:Y\ra J^1Y$ of the first jet prolongation
$J^1Y\ra Y.$
\end{defin}
\noindent
Further generalization of this idea leads us to the definition of $r$-th order connection, which is a section of $r$-th order jet prolongation of
a fibered manifold. According to the character of the target space we distinguish holonomic, semiholonomic and nonholonomic general connections. The coordinate form of a second order nonholonomic connection
$\Delta:Y\ra \wdt{J}^2Y$ is given by
$$y^p_i=F^p_i(x,y),\quad y^p_{0i}=G^p_i(x,y), \quad y^p_{ij}=H^p_{ij}(x,y),$$
where $F,G,H$ are arbitrary smooth functions. In case of linear connections all functions are linear in fiber coordinates.

Let us now recall the geometric concept of a
connection and its extension to higher order connections. The
following section is based on the paper \cite{rahula}.
\begin{defin} A connection on bundle $\pi:M_1\rightarrow
M$ is defined by the structure
$\triangle_h\oplus\triangle_v$ on a manifold $M_1$ where $\triangle_v=\ker T\pi$ is {\it vertical distribution}
tangent to the fibers and $\triangle_h$ is {\it horizontal distribution} complementary to the distribution
 $\triangle_v$\,. The transport of the fibers along the path $\gamma\subset M$ is realized by the horizontal lifts given by the distribution
 $\triangle_h$ on the surface $\pi^{-1}(\gamma)$.
 If the bundle is a vector one and the transport of fibers along an arbitrary path is linear, then
 the connection is called linear.
 \end{defin}

 We will assume that the base manifold $M$ is of dimension $n$ and the fibers are of dimension $r$. Then
 $$
\dim\triangle_h=n\,,\quad \dim\triangle_v=r\,.
 $$\\

On the neighborhood $U\subset M_1,$ let us consider local base and fiber coordinates:  $$(u^i, u^\alpha)\,,\, i=1,2,\ldots, n\,;\,
\alpha=n+1,\ldots, n+r.$$  Base coordinates $(u^i) \,\,$
are determined by the projection $\pi$ and the coordinates $(\bar{u}^i)$ on a neighborhood
$\bar{U}=\pi(U),\, u^i=\bar{u}^i\circ\pi\,.$\\

\begin{defin} On a neighborhood  $U\subset M_1$
we define a local
$($adapted$)$ basis of the structure
$\triangle_h\oplus\triangle_v\,$,
\begin{equation*}\label{adaptedB}
 (X_i \; X_\alpha)=\left ( {\partial \over\partial u^j} \, {\partial \over \partial u^\beta} \right )
 \cdot
 \begin{pmatrix}
  \delta^j_i & 0 \\
  \Gamma^\beta_i & \delta^\beta_\alpha \\
 \end{pmatrix},\quad
 \begin{pmatrix}
  \omega^i \\
  \omega^\alpha \\
 \end{pmatrix}=
 \begin{pmatrix}
 \delta^i_j & 0\\
 - \Gamma^\alpha_j & \delta^\alpha_\beta
 \end{pmatrix}\cdot
\begin{pmatrix}
 du^{j}\\
 du^{\beta}
 \end{pmatrix}.
\end{equation*}
The horizontal distribution $\triangle_h\,$ is the linear span of
the vector fields $(X_i)$ and the annihilator of the forms
$(\omega^\alpha)$,
\[
X_i=\partial_i+\Gamma^\beta_i\partial_\beta,\quad
\omega^\alpha=du^\alpha-\Gamma^\alpha_idu^i.
\]
\end{defin}

 \begin{defin}\label{class}A classical affine connection on manifold $M$
 is seen as a linear connection on the bundle $\pi_1:TM\rightarrow M$.
 On the tangent bundle $TM\rightarrow M$ one can define the structure $\triangle_h\oplus\triangle_v$.  The indices in the formulas are denoted by Latin letters
 all of them ranging from 1 to $n$. The functions $\Gamma^\alpha_i\,,
X_i\,, \omega^\alpha$
 are of the form $($in $\Gamma^\alpha_i$ the sign is changed to comply with the classical theory$)$\,$:$
 \begin{align*}
\Gamma^\alpha_i\quad&\rightsquigarrow\quad
-\Gamma^i_{jk}u^k_1\,,\\
X_i=\partial_i+\Gamma^\alpha_i\partial_\alpha\quad&\rightsquigarrow\quad
X_i=\partial_i-\Gamma^k_{ij}u^i_1\partial^1_k\,,\\
\omega^\alpha=du^\alpha-\Gamma^\alpha_idu^i\quad&\rightsquigarrow\quad
U^i_{12}=u^i_{12}+\Gamma^i_{jk}u^k_1u^j_2\,.
\end{align*}
\end{defin}

\begin{defin}\label{higher order} Higher order connections are defined as follows:
 on tangent bundle $TM$  the structure $\triangle\oplus\triangle_1$ is defined where $\ker T\rho_1=\triangle_1$,
 on $T(TM)$  the structure $\Delta\oplus\Delta_1\oplus\Delta_2\oplus\Delta_{12}$ is defined where
 $\ker T\rho_s=\Delta_s\oplus\Delta_{12}\,,
 s=1,2,$ etc.
 \end{defin}

\section{Connections on two-fold fibered manifolds}
More generally, one can define a second order connection by means of a two-fold fibered manifold. Note that the Definition \ref{higher order} is a special case of the following. A two-fold fibered manifold is a commutative diagram
\footnotesize{\[ \mathcal{M}
\]
\[
^{\rho_2}\swarrow\quad\qquad  \quad\searrow\,\,^{\rho_1}
\]
\begin{align}
\mathcal{M}_1 \qquad\qquad\qquad\qquad \mathcal{M}_2
\end{align}
\[
_{\pi_1}\searrow \qquad\qquad \swarrow\,\,_{\pi_2}
\]
\[
 M
\]}\normalsize where $\rho_1,\rho_2$ and $\pi_1,\pi_2$ -- four fibered manifolds
\[
 \dim
M=n,\,\,\,\dim\mathcal{M}_1=n+r_1,\,\,\,\dim\mathcal{M}_2=n+r_2,\,\,\,\dim\mathcal{M}=n+r_1+r_2+r_{12}.
\]
The double projection
$$\pi=\pi_1\circ\rho_2=\pi_2\circ\rho_1 : \mathcal{M}\rightarrow M$$
divides a manifold $\mathcal{M}$ to $n$-parameter family
of fibers of dimensions $(r_1+r_2+r_{12})$. Each fiber carries structure of another two fibers of dimensions $r_1+r_{12}$\, and \,$r_2+r_{12}$\, and these two fibers have the common intersection
of dimension $r_{12}$.

A two-fold fibered manifold is called a vector bundle if both fibrations
 $\pi_1,$ $ \pi_2, \rho_1$ è $\rho_2$ -- form vector bundles.

An example of a two-fold fibered manifold is the second order tangent bundle $T^2M$ of a manifold $M$. In this case $n=r_1=r_2=r_{12}$.\\

\begin{defin}
 A connection on a two-fold fibered manifold is defined by a structure on a manifold
$\mathcal{M}$:
\begin{equation}\label{structure}
\Delta\otimes\Delta_1\otimes\Delta_2\otimes\Delta_{12}\,,
\end{equation}
\[
\dim\Delta=n,\quad
\dim\Delta_1=r_1\,,\quad\dim\Delta_2=r_2\,,\quad\dim\Delta_{12}=r_{12}\,,
\]
$${\rm Ker} T\rho_2=\Delta_2\oplus\Delta_{12}\,,\qquad{\rm Ker} T\rho_1=\Delta_1\oplus\Delta_{12}$$
$$T\rho_2(\Delta\oplus\Delta_1)=T\mathcal{M}_1,\quad T\rho_1(\Delta\oplus\Delta_2)=T\mathcal{M}_2\,,$$
$$T\pi\Delta=TM.$$
\end{defin}

\begin{remark}
A connection on a two-fold vector fibered manifold is called linear if the structure \eqref{structure} induces on the manifolds $\pi_1, \pi_2, \rho_1$
and $\rho_2$ linear connections.
\end{remark}

\begin{remark}
Similarly, one can define a connection on a $k$--fold fibered manifold. In such case the commutative diagram would be represented by a $k$--dimensional cube. These manifolds would correspond to the $k$--th tangent bundle $T^kM$ of a manifold  $M$.
\end{remark}

On the neighborhoods
$$\mathcal{U}\subset\,\mathcal{M},\,\,\,
\mathcal{U}_1=\rho_2(\mathcal{U})\subset\,\mathcal{M}_1,\,\,\,
\mathcal{U}_2=\rho_1(\mathcal{U})\subset\,\mathcal{M}_2,\,\,\,
U=\pi(\mathcal{U})\subset\,M$$   we have the coordinate systems\,\,\,
$(u^i,\,\, u^{\alpha_1}, \,\, u^{\alpha_2},\,\,
u^{\alpha_{12}}),\,\,\,(u^i,\,\, u^{\alpha_1}),\,\,\,(u^i,\,\,
u^{\alpha_1}),\,\,\,(u^i).$\\

The transformation of coordinates on the neighborhoods
$\mathcal{U}$,
\begin{equation}(u^i,\, u^{\alpha_1},\,
 u^{\alpha_2},\, u^{\alpha_{12}})\,\rightsquigarrow\,(\tilde{u}^i,\, \tilde{u}^{\alpha_1},\,
 \tilde{u}^{\alpha_2},\, \tilde{u}^{\alpha_{12}})=(a^i,\, a^{\alpha_1},\,
 a^{\alpha_2},\, a^{\alpha_{12}}),\end{equation}
gives a Jacobi matrix:
\begin{equation}\left(
  \begin{array}{cccc}
    a^i_j & 0 & 0 & 0 \\
    a^{\alpha_1}_j & a^{\alpha_1}_{\beta_1} & 0 & 0 \\
    a^{\alpha_2}_j &0& a^{\alpha_2}_{\beta_2} & 0 \\
    a^{\alpha_{12}}_j & a^{\alpha_{12}}_{\beta_1}&a^{\alpha_{12}}_{\beta_2} & a^{\alpha_{12}}_{\beta_{12}} \\
  \end{array}
  \right).\end{equation}

See \cite{R1},\cite{R3}. Let us mention that the local (adapted) basis of such decomposition is represented by a matrix of the form
\begin{equation}\label{matrconHO}\left(
  \begin{array}{cccc}
    \delta^i_j & 0 & 0 & 0 \\
    \Gamma^{\alpha_1}_j & \delta^{\alpha_1}_{\beta_1} & 0 & 0 \\
    \Gamma^{\alpha_2}_j &0& \delta^{\alpha_2}_{\beta_2} & 0 \\
    \Gamma^{\alpha_{12}}_j & \Gamma^{\alpha_{12}}_{\beta_1}&\Gamma^{\alpha_{12}}_{\beta_2} & \delta^{\alpha_{12}}_{\beta_{12}} \\
  \end{array}
\right).\end{equation} The dual reper is given by the system of 1--forms:

\begin{align*}
\omega^i&=du^i,\nonumber\\
\omega^{\alpha_1}&=du^{\alpha_1}-\Gamma^{\alpha_1}_idu^i,\nonumber\\
\omega^{\alpha_2}&=du^{\alpha_2}-\Gamma^{\alpha_2}_idu^i,\\
\omega^{\alpha_{12}}&=du^{\alpha_{12}}-\Gamma^{\alpha_{12}}_{\alpha_1}
du^{\alpha_1}-\Gamma^{\alpha_{12}}_{\alpha_2}
du^{\alpha_2}-\bar{\Gamma}^{\alpha_{12}}_idu^i,\nonumber\\
\text{where}\quad
&\Gamma^{\alpha_{12}}_i-\bar{\Gamma}^{\alpha_{12}}_i=\Gamma^{\alpha_{12}}_{\beta_1}\,\Gamma^{\beta_1}_i
+\Gamma^{\alpha_{12}}_{\beta_2}\,\Gamma^{\beta_2}_i.
\end{align*}

\noindent
In case of linear connection the elements of the matrix \eqref{matrconHO} are of the form
\begin{align*}
\Gamma^{\alpha_1}_j=\Gamma^{\alpha_1}_{j\beta_1}u^{\beta_1}&,\quad\Gamma^{\alpha_2}_j=
\Gamma^{\alpha_2}_{j\beta_2}u^{\beta_2},\\
\Gamma_{\beta_1}^{\alpha_{12}}=\Gamma^{\alpha_{12}}_{{\beta_1}{\beta_2}}u^{\beta_2},
&\quad\,\,\,\Gamma_{\beta_2}^{\alpha_{12}}=\Gamma^{\alpha_{12}}_{{\beta_2}{\beta_1}}u^{\beta_1},\\
 \Gamma^{\alpha_{12}}_j=\Gamma^{\alpha_{12}}_{j{\beta_1}{\beta_2}}u^{\beta_1}
u^{\beta_2}+\Gamma^{\alpha_{12}}_{j{\beta_{12}}}u^{\beta_{12}}&,\quad
 \bar{\Gamma}^{\alpha_{12}}_j=\bar{\Gamma}^{\alpha_{12}}_{j{\beta_1}{\beta_2}}u^{\beta_1}
u^{\beta_2}+\bar{\Gamma}^{\alpha_{12}}_{j{\beta_{12}}}u^{\beta_{12}},\end{align*}
$$\Gamma^{\alpha_{12}}_{j{\beta_1}{\beta_2}}-\bar{\Gamma}^{\alpha_{12}}_{j{\beta_1}{\beta_2}}=
\Gamma^{\alpha_{12}}_{{\gamma_2}{\beta_1}}\Gamma^{\gamma_2}_{j\beta_2}
,$$ where the coefficients depend on the base coordinates $u^i$ only.

\section{Ehresmann prolongation}\label{ehres}

We show that Ehresmann prolongation plays an important role in determining all natural operators transforming
first order connections into higher order connections. Let us note that also natural transformations of semiholonomic
jet prolongation functor $\overline J^r$ are involved. To find the details about this topic we refer to \cite{kms},\cite{kolmod}, \cite{vasik}.
For our purposes, it is enough to consider $r=2.$ We use the notation of \cite{kms}, where the map $e:\overline{J}^2Y\ra \overline{J}^2Y$ is obtained from the natural exchange map
$e_{\Lambda}:J^1J^1Y\ra J^1J^1Y$ as a restriction to the subbundle $\overline{J}^2Y\subset J^1J^1Y$.
Note that while $e_{\Lambda}$ depends on the linear connection $\Lambda$ on $M$, its restriction $e$ is independent of
any auxiliary connections. We remark, that originally the map $e_{\Lambda}$ was introduced by M. Modugno. We also recall that
J. Pradines introduced a natural map $\overline{J}^2Y\ra \overline{J}^2Y$ with the same coordinate
expression.

Now we are ready to recall the following assertion, see \cite{vasik} for the proof.

\begin{prop}
\label{ehres}
All natural operators transforming first order connection $\Gamma:Y\ra J^1 Y$ into second
order semiholonomic connection $Y\ra\overline{J}^2 Y$ form a one parameter family
\begin{equation*}
\Gamma\mapsto k\cdot(\Gamma\ast\Gamma)+(1-k)\cdot e(\Gamma\ast\Gamma),\qquad k\in\mathbb{R}.
\end{equation*}
\end{prop}

This shows the importance of Ehresmann prolongation in the theory of prolongations of connections.
Let us now recall a more general concept of a product of two connections:

Given two
higher order connections $\Gamma:Y\rightarrow\widetilde{J}^{r}Y$
and $\overline{\Gamma}:Y\rightarrow \widetilde{J}^{s}Y,$ the
product of $\Gamma$ and $\overline{\Gamma}$ is the $(r+s)$-th
order connection $\Gamma\ast\overline{\Gamma}:Y\rightarrow
\widetilde{J}^{r+s}Y$ defined by
$$\Gamma\ast\overline{\Gamma}=\widetilde{J}^{s}\Gamma\circ\overline{\Gamma}.$$

Particularly, if both $\Gamma$ and $\overline{\Gamma}$ are of the first order,
then $\Gamma\ast\overline{\Gamma}:Y\rightarrow \widetilde{J}^{2}Y$
is semiholonomic if and only if $\Gamma=\overline{\Gamma}$ and
$\Gamma\ast\overline{\Gamma}$ is holonomic if and only if $\Gamma$
is curvature-free, \cite{kol2}, \cite{virsik}.

As an example we show the coordinate expression of an arbitrary nonholonomic second order connection
and of the product of two first order connections. The coordinate form of $\Delta:Y\ra \wdt{J}^2Y$ is
$$y^p_i=F^p_i(x,y),\quad y^p_{0i}=G^p_i(x,y), \quad y^p_{ij}=H^p_{ij}(x,y),$$
where $F,G,H$ are arbitrary smooth functions. Further, if the coordinate expressions of two first order connections
$\Gamma,\overline{\Gamma}:Y\ra J^1Y$ are
\begin{align}\label{connclas}\Gamma:\quad y^p_i=F^p_i(x,y),\qquad \overline{\Gamma}:\quad y^p_i=G^p_i(x,y),\end{align}
then the second order connection $\Gamma\ast\overline{\Gamma}:Y\ra\wdt{J}^2Y$ has equations
\begin{align}\label{equation}
y^p_i = F^p_i, \quad
y^p_{0i} = G^p_i,\quad
y^p_{ij} = \frac{\partial F^p_i}{\partial x^j}+\frac{\partial F^p_i}{\partial y^q}G^q_j.
\end{align}
For linear connections, the coordinate form would be obtained by substitution
\begin{align*}
F^p_i&=F^p_{iq}y^q,\\
G^p_i&=G^p_{iq}y^q
\end{align*}
in the equations \eqref{connclas}, where $F^p_{iq}$ and $G^p_{iq}$
are functions of the base manifold coordinates $x_i.$ For order
three see \cite{vasik1}.\\

\section{Tangent funktor and Ehresmann prolongation}

If we apply the tangent functor $T$ two times on a projection $\pi:E\rightarrow M$ and a section $\sigma:
M\rightarrow E$ we obtain
$$T\pi:TE\rightarrow TM\,,\,\,T^2\pi:T^2E\rightarrow T^2M,$$
$$T\sigma: TM\rightarrow
TE,\,\,T^2\sigma: T^2M\rightarrow T^2E,$$
respectively.
The mappings $\sigma,
T\sigma$ and $T^2\sigma$ define the sections of fibered manifolds $\pi, T\pi$ and $
T^2\pi.$

Let us consider local coordinates on the following manifolds in the form $$\text{on}\quad M,\, TM,\,
T^2M\,:\,\,\,(x^i),\,\,(x^i, x^i_1),\,\, (x^i, x^i_1, x^i_2,
x^i_{12}),$$
$$\text{and on}\quad E,\, TE,\,
T^2E\,:\,\,\,(y^p),\,\, (y^p, y^p_1),\,\, (y^p, y^p_1, y^p_2,
y^p_{12}).$$ Let us also consider for a function $f$ defined on a manifold $M$, its following differentials on $T^2M$ in local coordinate form:
\[
f_1\doteq f_ix^i_1,\,\,\,f_2\doteq f_ix^i_2,\,\,\,f_{12}\doteq
f_{ij}x^i_1x^j_2+f_ix^i_{12},\,\,\text{where}\,\, f_i=\frac{\partial
f}{\partial x^i}\,,\,\,f_{ij}=\frac{\partial^2 f}{\partial
x^i\partial x^j}.
\]
Furthermore,
$f_1=df\circ\rho_1\,,\,\,f_2=df\circ\rho_2\,,\,\,f_{12}=d\,^2f$.
We use these notations in the formulae bellow.

If the section $\sigma$ is defined by local functions $\Gamma^p$, then
the sections $T\sigma$ and $T^2\sigma$ are defined by its differentials
$\Gamma_1^p$\,,\,\,$\Gamma_2^p$\,\,\,and\,\,\,$\Gamma_{12}^p,$
\begin{align}
&\sigma:\,\,x^i\rightsquigarrow y^p=\Gamma^p,\nonumber\\
& T\sigma:\,\, (x^i,x^i_1)\rightsquigarrow
(y^p,y^p_1)=(\Gamma^p,\Gamma^p_1),\nonumber\\
&T^2\sigma:\,\,(x^i, x^i_1, x^i_2, x^i_{12})\rightsquigarrow (y^p,
y^p_1, y^p_2, y^p_{12})=(\Gamma^p, \Gamma^p_1, \Gamma^p_2,
\Gamma^p_{12}),\nonumber\\
& \label{coords}\text{where}\quad\Gamma^p_1=\Gamma^p_ix^i_1,\quad\Gamma^p_2=\Gamma^p_ix^i_2,
\quad\Gamma^p_{12}=\Gamma^p_{ij}x^i_1x^j_2+\Gamma^p_ix^i_{12}.
\end{align}

The case when the coefficients $\Gamma^p_i,\,\,\Gamma^p_{ij}$ in \eqref{coords} are arbitrary functions, corresponds to a nonholonomic connection on the fibered manifold $\pi$.

The case when $\Gamma^p_{ij}=\displaystyle
{\frac{\partial \Gamma^p_i}{\partial x^j}}$,\, where\,\, $\Gamma^p_i$
are arbitrary functions corresponds to a semiholonomic connection on the fibered manifold $\pi$.

The case when
$\Gamma^p_1=d\Gamma^p\circ\rho_1\,,\,\,\Gamma^p_2=d\Gamma^p\circ\rho_2\,,\,\,\Gamma^p_{12}=d\,^2\Gamma^p\,,\,\,$
corresponds to a holonomic connection on the fibered manifold $\pi$.

The functions $\Gamma^p_i,\,\Gamma^p_{ij}$ define nonholonomic, semiholonomic or holonomic Ehresmann prolongation of a connection, respectively.

\begin{remark}
Nonholonomic prolongation induces a connection on a double fibered manifold
$$
J\,\rightarrow\, E\,\rightarrow\,
M\,:\,\,\,y^p_i\,\rightsquigarrow\,y^p\,\rightsquigarrow\, x^i.
$$
On the fibered manifold $E\,\rightarrow\,M$ the fiber transformations are given by the Pfaff system
\[
\omega^p\equiv dy^p - \Gamma^p_idx^i=0,
\]
more precisely, along a curve $x^i(t)$ -- by the system of first order ODEs
\begin{align}\label{ode} \dot{y}^p=\Gamma^p_i\,\dot{x}^i.
\end{align}
In case $(\Gamma^p_{12}, x^i_1, x^j_2,
x^i_{12})\,\rightsquigarrow\,(\ddot{y}^p\,,\,\dot{x}^i\,,\,\dot{x}^j\,,\,\ddot{x}^i)$
we obtain the system of second order ODEs:
\[
\Gamma^p_{12}=\Gamma^p_{ij}x^i_1x^j_2+\Gamma^p_ix^i_{12}\,\,\rightsquigarrow\,\,
\ddot{y}^p=\Gamma^p_{ij}\dot{x}^i\dot{x}^j+\Gamma^p_i\ddot{x}^i.
\]
Considering the system \eqref{ode}, we obtain for fiber coordinates $y^\alpha,
y^\alpha_i$ system of first order ODEs \begin{align}
\left\{\begin{array}{ll}
    \dot{y}^p &=\,\Gamma^p_i\,\dot{x}^i,\\
    \dot{y}^p_i &=\,\Gamma^p_{ij}\,\dot{x}^j.
\end{array} \right.\end{align}
The sections of fibers along a curve $x^i(t)$ are given.

\noindent The horizontal distribution $\triangle_h$ is $n$-dimensional and described by the vector field
\[
X_i=\partial_i+\Gamma^p_i\,\partial_p+\Gamma^p_{ij}\,\partial^j_p\,,\quad
\text{where}\quad\partial_i=\frac{\partial}{\partial x^i},\,\,
\partial_p=\frac{\partial}{\partial y^\alpha},\,\,\partial^j_p=\frac{\partial}{\partial
y^p_j}\,.
\]\\

\end{remark}

\vspace{1cm}
\noindent
MAIDO RAHULA\\
University of Tartu,\\
Institute of Mathematics\\
J.Liivi 2\\
50409 Tartu, Estonia\\
\noindent e-mail:{\tt \ rahula@ut.ee} \vspace{0.5cm}
\newline
PETR VA\v S\'IK\\
Institute of Mathematics\\
Brno University of Technology\\
FME BUT Brno, Technick\'a 2\\
616 69 Brno, Czech Republic\\
\noindent e-mail:{\tt \ vasik@fme.vutbr.cz }


\begin{thebibliography}{2}
\bibitem{R1}  Atanasiu, G., Balan, V., Br\^{\i}nzei, N., Rahula, M.:
{\it Differential Geometric Structures: Tangent Bundles, Connections
in Bundles, Exponential Law in the Jet Space}  (in Russian),
Librokom, Moscow, 2010.

\bibitem{ehres} Ehresmann Ch.: {\it Les connexions infinit\'{e}simales dans un
        espace fibr\'{e} dif\-f\'{e}rentiable}, Coll. de Topologie, Bruxelles,
        CBRM (1950), 29-55;\,\, {\OE}uvres Compl\`{e}tes, t.\,I, 28, 179-204.

\bibitem{ehres1} Ehresmann Ch., {\it Cat\'{e}gories doubles et cat\'{e}gories
structur\'{e}es}, C.R. Acad. Sci., Paris, 256(1958), 1198-1201.

\bibitem{e} Epstein, M.:
\emph{The Geometrical Language of Continuum Mechanics},
Cambridge University Press, 2010.

\bibitem{kol2} Kol\'a\v r I.: \textit{On the torsion of spaces with connections}, Czechoslovak Math. J. 21(1971), 124-136.
\bibitem{kms} Kol\'a\v r, I., Michor P. W., Slov\'ak J.: \textit{Natural Operations in Differential Geometry,} Springer-Verlag, 1993.
\bibitem{kolmod} Kol\'a\v r, I., Modugno, M.: \textit{Natural maps on the iterated jet prolongation of a fibered manifold}, Annali
di Matematica CLVIII (1991), 151-165.
\bibitem{mik1} Mikulski, W. M.: \emph{Some natural operations on vector fields}, Rend. Math. (Serie VII) 12(1992), 783-803.
\bibitem{pradines} Pradines J.: {\it Suites  exactes vectorielles doubles et
connexions}, C.R. Acad. Sci., Paris, 278(1974), 1587-1590.

\bibitem{R3} Rahula, M.: {\it New Problems in Differential Geometry}, WSP,
1993.

\bibitem{rahula} Rahula, M., Va\v s\'ik, P., Voicu, N.: \emph{Tangent structures: sector-forms, jets and connections}, J. Phys.: Conf. Ser. 346 (2012) 012023.

\bibitem{vasik} Va\v s\'ik, P.: \textit{On the Ehresmann Prolongation,}
Annales Universitatis Mariae Curie Sklodowska, LXI (2007), Sectio A, 145-153.

\bibitem{vasik1} Va\v s\'{i}k, P.: \textit{Transformations of semiholonomic 2- and 3-jets and semiholonomic
prolongation of connections,} Proc. Est. Acad. Sci. 59 (2010), 4, 375-380.

\bibitem{virsik} Virsik G.: \textit{On the holonomity of higher order connections}, Cahiers Topol. G\'eom. Diff.
12(1971), 197-212.

\end{thebibliography}
\end{document}